\documentclass[11pt]{article}

\usepackage{amsmath}
\usepackage{amssymb}
\usepackage{eucal}

\begin{document}

\baselineskip 16pt

\title{Some new characterizations of $PST$-groups}

\author{Xiaolan Yi \\
{\small Department of Mathematics, Zhejiang Sci-Tech University,}\\
 {\small Hangzhou 310018, P.R.China}\\
{\small E-mail:
yixiaolan2005@126.com}\\ \\
{ Alexander N. Skiba }\\
{\small Department of Mathematics, Francisk Skorina Gomel State University,}\\
{\small Gomel 246019, Belarus}\\
{\small E-mail: alexander.skiba49@gmail.com}}

\date{}
\maketitle

\begin{abstract}   Let  $H$ and $B$ be  subgroups  of a finite group 
  $G$ such that $G=N_{G}(H)B$.   
 Then we say that   $H$ is \emph{quasipermutable} (respectively 
\emph{$S$-quasipermutable}) in $G$ 
 provided  $H$  permutes  with $B$ and
with every subgroup (respectively  with every Sylow  subgroup) 
 $A$ of $B$  such that $(|H|, |A|)=1$.  
 In this paper we  analyze the influence of 
$S$-quasipermutable  and quasipermutable 
 subgroups on the structure of $G$.   As an application, we give new 
characterizations of soluble $PST$-groups.

\end{abstract}

\let\thefootnoteorig\thefootnote
\renewcommand{\thefootnote}{\empty}


\let\thefootnote\thefootnoteorig

\section{Introduction}

Throughout this paper, all groups are finite and $G$ always denotes
a finite group.   Moreover $p$ is always supposed to be a prime and
$\pi$  is a  subset of the set $\Bbb{P}$ of all primes; $\pi (G)$ denotes the
 set of all primes dividing $|G|$.

A    subgroup $H$ of $G$ is said  to be \emph{quasinormal}  or
\emph{permutable}  in $G$ if $H$ permutes with  every subgroup $A$ of $G$, that is, $HA=AH$; 
 $H$ is          said to be \emph{$S$-permutable} 
 in $G$ if $H$ permutes with  every Sylow subgroup of $G$.

A group   $G$ is called a     \emph{$PT$-group}  if
 permutability   is a 
transitive relation  on $G$, that is, every  permutable subgroup of
 a permutable  subgroup 
 of $G$ is permutable  in $G$.  
A group   $G$ is called a     \emph{$PST$-group} if
 $S$-permutability is a         transitive relation  on $G$.

As well as  $T$-groups, $PT$-groups and $PST$-groups possess many interesting 
properties (see  Chapter 2 in  \cite{prod}).
 The  general description of
$PT$-groups and $PST$-groups  were firstly obtained by Zacher 
\cite{G.Zacher} and  Agrawal \cite{Agr},
  for the soluble  case, and 
  by Robinson in \cite{217}, for the general case. Nevertheless,
in the  further publications, the  authors (see   for example recent papers\cite{78}-- 
\cite{khaledII}) have found out and  described 
many other   interesting characterizations   of soluble $PT$ and $PST$-groups.

 In this paper we  give new   "Hall"-characterizations of soluble 
$PST$-groups on the basis of  the following

{\bf Definition 1.1.}  We   say that  
  a subgroup     $H$ is \emph{quasipermutable} (respectively 
\emph{$S$-quasipermutable}) in $G$ 
 provided  $H$  permutes  with $B$ and
with every subgroup (respectively  with every Sylow  subgroup) 
 $A$ of $B$  such that $(|H|, |A|)=1$.

Examples and some applications of quasipermutable 
 subgroups were discussed in our  
 papers  \cite{Bull} and  \cite{proble} (see also remarks in Section 5 below).
In this paper,
 we give  the following result, which we consider as one more
  motivation for introducing the concept of quasipermutability.

{\bf Theorem  A.}   {\sl Let  $D=G^{\cal N} $ and   $\pi =\pi (D)$. 
 Then the following 
statements are equivalent:}   

(i)  {\sl $D$  is a Hall subgroup of $G$ and 
every Hall subgroup of $G$ is   quasipermutable in $G$.}

(ii) {\sl $G$ is a soluble $PST$-group.}

(iii)  {\sl Every subgroup of $G$ is  quasipermutable in $G$.}

(iv)  {\sl Every $\pi$-subgroup of $G$ 
  and some minimal supplement of $D$ in 
$G$  are  quasipermutable in $G$.}

In the  proof    Theorem A 
we use  the next  three our results.

A subgroup $S$ of $G$ is called a \emph{Gasch\"utz} subgroup 
of $G$ (L.A. Shemetkov   \cite[IV,  15.3]{26}) if  $S$ is  supersoluble
 and for any subgroups $K \leq  H$ of $G$,  where $S\leq K$, the number $|H:K|$ 
is not prime.

{\bf Theorem B.} {\sl  The following statements are equivalent:}

(I) {\sl   $G$ is soluble, and  if  $S$ is a    Gasch\"utz     
 subgroup  of $G$, then 
every Hall subgroup $H$ of $G$  satisfying $\pi (H)\subseteq \pi (S)$ 
 is quasipermutable in $G$.}

(II)    {\sl  $G$ is supersoluble and the following hold: }

 (a) {\sl $G=DC$,  where $D=G^{\cal N}$
 is an abelian complemented subgroup of $G$ and  $C$ is a Carter subgroup  of 
$G$;}

(b) {\sl   $D\cap C$ is normal in $G$  and
 $(p, |D/D\cap C|)=1$ for all prime divisors $p$ of $|G|$  satisfying $(p-1, |G|)=1$.}

(c) {\sl   For any non-empty set $\pi $ of primes,  every $\pi $-element of  any
  Carter subgroup  of $G$     induces a power      automorphism on the
 Hall $\pi'$-subgroup of $D$.}

(III) {\sl Every Hall subgroup of $G$ is quasipermutable in $G$.}

Let $\cal F$ be a class of groups. If $1\in {\cal F}$, then we write $G^{\cal F}$ to denote the
intersection of all normal subgroups $N$ of $G$ with  $G/N\in {\cal
F}$.   The class $\cal F$ is said to be a \emph{formation} if either
${\cal F}= \varnothing $ or $1\in {\cal F}$ and every
homomorphic image of $G/G^{\cal F}$  belongs to $ {\cal F}$ for any group $G$. 
 The formation ${\cal F}$
is said to be  \emph{saturated} if $G\in {\cal F}$ whenever $G/\Phi (G) \in {\cal 
F}$.    A subgroup $H$ of $G$ is said to be an \emph{$\cal F$-projector}  of $G$ 
provided $H\in {\cal F}$  and $E=E^{\cal F}H$ for any subgroup $E$ of $G$ 
containing $H$. By the Gasch\"utz's theorem \cite[VI, 9.5.4 and 9.5.6]{hupp}, for any   saturated formation
$\cal F$,  every  soluble group $G$ has  an  $\cal F$-projector    and any two
 $\cal F$-projectors  of $G$ are conjugate.

{\bf Theorem  C.}   {\sl Let $\cal F$ be a  saturated formation containing 
all nilpotent groups. Suppose that $G$ is soluble and  let $\pi =\pi (C) 
\cap  \pi( G^{\cal F})$, where $C$ is an $\cal F$-projector of $G$. If 
every  maximal subgroup of every Sylow $p$-subgroup of $G$ is 
 $S$-quasipermutable in  $G$ for all  $p\in \pi $, then  $G^{\cal F}$ is a 
Hall subgroup of $G$.}

{\bf Theorem  D.}   {\sl  Let $\cal F$ be a saturated formation containing 
all supersoluble  groups and   $\pi =\pi (F^{*}(G^{\cal F}))$.  If 
$G^{\cal F}\ne 1$, then for some  $p\in \pi $  some maximal subgroup  of a Sylow
 $p$-subgroup of $G$  is not   $S$-quasipermutable in  $G$.}

In this theorem   $F^{*}(G^{\cal F})$ denotes the generalized  Fitting 
subgroup of  $G^{\cal F}$, that is, the product of all normal quasinilpotent subgroups of  
 $G^{\cal F}$. 

The main tool in the proofs of Theorems C and D is the following  our result.

{\bf Proposition.}   {\sl Let $E$   be a normal subgroup of $G$  and
 $P$  a Sylow $p$-subgroup of $E$ such that $|P| > p$. }

(i) {\sl  If  every
 number $V$ of some fixed  ${\cal M}_{\phi}(P)$ is $S$-quasipermutable in
$G$, then $E$ is $p$-supersoluble.  }

(ii) {\sl  If 
 every maximal subgroup of $P$  is $S$-quasipermutable in
$G$, then every     chief factor of $G$   between $E$ and $O_{p'}(E)$ 
  is cyclic.   }

(iii) {\sl If every maximal subgroup of every  Sylow subgroup of $E$ is $S$-quasipermutable 
in $G$, then  every chief factor of $G$   below  $E$ is cyclic. }

In this proposition   we write ${\cal M}_{\phi}(G)$, by analogy with \cite{shirong}, 
 to denote a set of      maximal subgroups of $G$
such that  ${\Phi}(G)$ coincides with the intersection of all
subgroups in ${\cal M}_{\phi}(G)$.

Note that Proposition may be  independently interesting because this result 
 unifies and generalize  
 many known results, and in particular, Theorems 1.1--1.5 in  
 \cite{shirong} (see Section 5).    In Section 5 we discus also some
 further applications of the results.

All unexplained notation and terminology are standard. The reader is
referred to \cite{26},  \cite{DH},    or  \cite{Bal-Ez} if necessary.

\section{Basic    Propositions}

  Let  $H$ be a subgroup  of  $G$.
 Then we say, following \cite{Bull},  that  $H$ is    \emph{propermutable}
   (respectively \emph{$S$-propermutable}) 
  in $G$ provided there is  a subgroup $B$ of
$G$ such that $G=N_{G}(H)B$      and  $H$  permutes with all  subgroups (respectively  with all
 Sylow  subgroups)   of $B$.

  {\bf Proposition  2.1.} {\sl  Let  $H\leq G$  and $N$ a normal subgroup  of $G$. Suppose that
$H$ is  quasipermutable ($S$-quasipermutable) in $G$.}

(1) {\sl If either $H$ is a Hall subgroup of $G$ or 
  for every prime $p$ dividing $|H|$ and for every Sylow $p$-subgroup $H_{p}$
of $H$ we have $H_{p}\nleq N$, then $HN/N$ is quasipermutable
 ($S$-quasipermutable, respectively)  in $G/N$.    }

(2) {\sl If $\pi =\pi (H)$ and  $G$ is $\pi$-soluble, then $H$
 permutes with some Hall $\pi'$-subgroup of $G$.  }

(3) {\sl   $H$ permutes with some Sylow $p$-subgroup
 of $G$ for every prime  $p$ dividing  $|G|$ such that   $(p, |H| )=1$.}

(4) {\sl  $|G:N_{G}(H\cap N)|$ is a
$\pi$-number,  where $\pi = \pi (N)\cup  \pi (H)$.}

(5) {\sl If $H$ is  propermutable ($S$-propermutable) in $G$, then 
  $HN/N$ is propermutable ($S$-propermutable, respectively) in $G/N$.    }

(6) {\sl If    $H$ is $S$-propermutable in $G$, then $H$ 
permutes with some Sylow $p$-subgroup
 of $G$ for any prime $p$ dividing $|G|$.   }

(7)  {\sl  Suppose that $G$ is $\pi$-soluble.  If $H$ is a Hall $\pi$-subgroup
 of $G$, then  $H$ is 
propermutable  ($S$-propermutable, respectively) in $G$.   }

{\bf Proof.}  By hypothesis,   there  is a subgroup $B$ of $G$ such that $G=N_{G}(H)B$
 and   $H$ permutes with  $B$ and with all 
   subgroups  (with all Sylow    subgroups, respectively) $A$ of  $B$
 such that  $(|H|, |A|)=1$.

(1)  It is clear that 
 $$G/N=(N_{G}(H)N/N)(BN/N)=N_{G/N}(HN/N)(BN/N).$$
 Let  $K/N$ be any  subgroup (any Sylow subgroup, respectively) of
$BN/N$   such that $(|HN/N|, |K/N|)=1$.   
  Then $K=(K\cap  B)N$.  Let $B_{0}$ be  a minimal supplement of $K\cap B\cap N$ to $K\cap B$.
Then $K/N=(K\cap 
B)N/N =B_{0}(K\cap B\cap N)N/N=B_{0}N/N$  and   $K\cap B\cap N\cap B_{0}=N\cap B_{0}\leq \Phi (B_{0})$.
Therefore $\pi (K/N)=\pi (K\cap B/K\cap B\cap N)=\pi (B_{0})$, so  
$(|HN/N|, |B_{0}|)=1$.    Suppose that  some prime  $p\in \pi (B_{0})$ 
divides $|H|$, and let $H_{p}$ be  
 a Sylow $p$-subgroup of $H$. We shall show that  
 $H_{p}\nleq N$. In fact, we may suppose that 
 $H$ is a Hall subgroup of $G$. But in this case,
  $H_{p}$ is 
 a Sylow $p$-subgroup of $G$. Therefore, since $p\in \pi (B_{0})\subseteq \pi (G/N)$,  
$H_{p}\nleq N$.   Hence $p$ 
divides  $|HN/N|$, a contradiction.  Thus  $(|H|, |B_{0}|)=1$, so in the case, when $H$ is 
quasipermutable in $G$,   we have  $HB_{0}=B_{0}H$ and hence   $HN/N$
permutes with  $K/N=B_{0}N/N$.  Thus $HN/N$ is quasipermutable in 
$G/N$.

 Finally, suppose that $H$  is $S$-quasipermutable in $N$. In this case, 
 $B_{0}$   is a $p$-subgroup of $B$, so
for some Sylow $p$-subgroup $B_{p}$ of $B$ we have   $B_{0}\leq  B_{p}$ and $(|H|, p)=1$. 
Hence 
$K/N=B_{0}N/N\leq B_{p}N/N$, which implies that $K/N= B_{p}N/N$.
 But $H$ permutes with $B_{p}$ by hypothesis, so    $HN/N$
permutes with  $K/N$.
 Therefore $HN/N$ is $S$-quasipermutable in $G/N$.

(2)   By \cite[VI, 4.6]{hupp},   there are
 Hall   $\pi'$-subgroups $E_{1}$, $E_{2}$ and $E$ of $N_{G}(H)$,
$B$ and $G$, respectively,   such that $E=E_{1}E_{2}$.
 Then $H$ permutes with all Sylow subgroups of  $E_{2}$ by hypothesis, so

$$HE=H(E_{1}E_{2})=(HE_{1})E_{2}=(E_{1}H)E_{2}=$$ $$
 E_{1}(HE_{2})=E_{1}(E_{2}H)=(E_{1}E_{2})H=EH$$
by \cite[A,  1.6]{DH}.

(3) See the proof of (2).

(4) Let $p$ be a prime such that $p\not \in \pi $.   Then by (3),
there is a  Sylow $p$-subgroup $P$ of $G$ such that $HP=PH$ is a
subgroup of $G$. Hence $HP\cap N=H\cap N$ is a normal subgroup of
$HP$. Thus $p$ does not divide $|G:N_{G}(H\cap N)|$.

(5) See the proof of (1).

(6)   See the proof of (2).

(7)   Since $G$ is 
$\pi$-soluble, $B$ is $\pi$-soluble.  Hence by   \cite[VI, 1.7]{hupp},
 $B=B_{\pi}B_{\pi'}$  
where  $B_{\pi}$   is a Hall   $\pi$-subgroup of $B$ and  $B_{\pi'}$   is a 
Hall   $\pi'$-subgroup of $B$.   By \cite[VI, 4.6]{hupp},    there are
 Hall $\pi$-subgroups $N_{\pi}$, $B_{\pi}$ and $G_{\pi}$ of $N_{G}(H)$,
$B$ and $G$, respectively,   such that $G_{\pi}=N_{\pi}B_{\pi}$. But since 
$H\leq N_{\pi}$, 
 $N_{\pi}$ is a Hall   $\pi$-subgroup of $G$. Therefore   
$G_{\pi}=N_{\pi}B_{\pi}=N_{\pi}$, so  $B_{\pi}\leq N_{\pi}$. 
Hence 
 $G=N_{G}(H)B=N_{G}(H)B_{\pi}B_{\pi'}=N_{G}(H)B_{\pi'}$, so  $H$ is 
propermutable ($S$-propermutble, respectively) in $G$.

A group $G$ is said to be a \emph{$C_{\pi }$-group} provided $G$ has a Hall 
${\pi }$-subgroup   and any two Hall ${\pi }$-subgroups of $G$ are 
conjugate. 

On the basis of Proposition 2.1  the following two results are  proved.

{\bf Proposition  2.2.} {\sl Let    $H$ be 
 a  Hall $S$-quasipermutable subgroup of  $G$. If   $\pi = \pi (|G:H|)$, then  $G$
 is a $C_{\pi }$-group.}

{\bf Proposition 2.3.} {\sl Let  $E$ be a normal subgroup of $G$ and  $H$ 
 a  Hall $\pi$-subgroup of  $E$.     
 If $H$ is nilpotent and $S$-quasipermutable in $G$, then $E$ is $\pi$-soluble.}

\section{Groups with a  Hall quasipermutable subgroup }

A group $G$ is said to be \emph{$\pi$-separable} if every chief factor of 
$G$ is either a $\pi$-group or  a $\pi'$-group. Every $\pi$-separable 
group $G$ has a series $$1=P_{0}(G)\leq M_{0}(G) < P_{1}(G) <  M_{1}(G) 
< \ldots    <  P_{t}(G)\leq M_{t}(G)=G   $$ such that
 $$M_{i}(G)/P_{i}(G) =O_{\pi'}(G/P_{i}(G))$$ ($i=0, 1,   \ldots , t$)
and  $$P_{i+1}(G)/M_{i}(G)=  O_{\pi}(G/M_{i}(G))$$  
($i=1, \ldots , t$)

The number $t$ is called the \emph{$\pi$-length} of $G$ and denoted by 
$l_{\pi}(G)$ (see \cite[p. 249]{rob}).

One more result, which we use use in the proof of our main results,  is the 
following 

{\bf Theorem 3.1.} {\sl  Let $H$ be a
  Hall subgroup of $G$ and $\pi =\pi (H)$. Suppose that  $H$ is 
quasipermutable in $G$.   }

(I) {\sl If $ p > q  $   for all primes $p$ and $q$ such that $p\in \pi $  and 
$q$ divides $|G:N_{G}(H)|$, then $H$ is normal in $G$.}

(II)  {\sl If $H$ is supersoluble, then $G$ is $\pi$-soluble.}

(III) {\sl If $H$ is  $\pi$-separable, then the following fold:}

(i) {\sl $H'\leq O_{\pi}(G)$. If, in addition,
$N_{G}(H)$ is nilpotent, then  $G'\cap H  \leq O_{\pi}(G)$.}

(ii) {\sl   $l_{\pi}(G) \leq 2$ and $l_{\pi'}(G) \leq 2$.  }

(iii) {\sl If for some prime $p\in  \pi'$ a Hall $\pi'$-subgroup $E$ 
 of $G$
 is $p$-supersoluble, then $G$ is   $p$-supersoluble. }

Let $\cal M$  and $\cal H$ be non-empty  formations. Then the 
\emph{product}  ${\cal M} {\cal H}$ of these formations is   the  class  of all 
groups $G$ such that  $G^{\cal H}\in {\cal M}$.  It is well-known that such an operation 
 on the set  of all non-empty formations is  associative  (Gasch\"utz). The symbol
 ${\cal M}^{t}$ denotes  the product of $t$ copies of ${\cal M}$.

We shall need  following well-known 
 Gasch\"utz-Shemetkov's theorem \cite[Corollary 7.13]{100}.                                 

{\bf Lemma 3.2}.  {\sl   The product of any two non-empty
 saturated formations is also a saturated formation.}
 
In in the proof of Theorem 3.1 we use the following

{\bf Lemma 3.3.} {\sl   The class  $\cal F$ of all   $\pi$-separable groups $G$ 
with  $l_{\pi}(G) \leq t$ is a saturated formation.}

{\bf Proof.}  It is not difficult to show that   for any  non-empty
 set $\omega \subseteq \Bbb{P}$ the 
class ${\cal G}_{\omega}$ of all  $\omega$-groups is a saturated   
formation and that ${\cal F}=({\cal G}_{\pi'}{\cal G}_{\pi})^{t}{\cal G}_{\pi'}$. 
 Hence ${\cal F}$ is a    saturated formation by Lemma 3.2.

{\bf Lemma 3.4.} {\sl Suppose that $G$  is separable. If   Hall 
$\pi$-subgroups
 of $G$ are abelian, then $l_{\pi}(G) \leq 1$.}

{\bf Proof.} Suppose that this lemma  is false and let 
$G$ be a counterexample of minimal order.    Let $N$ be a minimal normal subgroup of $G$.  
 Since $G$   is 
$\pi$-separable,   $N$ is a $\pi$-group or a $\pi'$-group.  
It is clear that the hypothesis holds for $G/N$, so   $l_{\pi}(G/N) \leq 1$ by the choice of $G$.
By Lemma 3.3,  the class  of all $\pi$-soluble groups with $l_{\pi}(G) \leq 1$ is 
a saturated formation.   Therefore $N$ is a unique minimal normal subgroup 
of $G$, $N\nleq \Phi (G)$ and $N$ is not a  $\pi'$-group. Hence $N$ is a  $\pi$-group and 
$N=C_{G}(N)$ by \cite[A, 15.2]{DH}.  Therefore  $N\leq H$, where $H$ is a Hall 
$\pi$-subgroup of $G$.  But since   $H$ is abelian, $N=H$ is a Hall 
$\pi$-subgroup of $G$.  Hence  
$l_{\pi}(G) \leq 1$. 

A group $G$ is called \emph{$\pi$-closed}  provided $G$ has a normal 
Hall $\pi$-subgroup.

{\bf Lemma 3.5.} {\sl  Let $H$ be a Hall $\pi$-subgroup of $G$. If $G$ is 
$\pi$-separable and $H\leq Z(N_{G}(H))$, then $G$ is $\pi'$-closed.   }

{\bf Proof.}  Suppose that this lemma  is false and let 
$G$ be a counterexample of minimal order.  Then $G\ne H$.
     The class $\cal F$  of all  $\pi'$-closed groups coincides with the 
product ${\cal G}_{\pi'}{\cal G}_{\pi}$. Hence $\cal F$   is a saturated formation by
 Lemma 3.2.  
 Let $N$ be a minimal normal subgroup of $G$. Since $G$  is 
$\pi$-separable,   $N$ is a $\pi$-group or a $\pi'$-group. Moreover, $G$ 
is a  $C_{\pi}$-group   by    \cite[9.1.6]{rob}), so the  hypothesis holds for $G/N$.
 Hence $G/N$ is 
$\pi'$-closed by the choice of $G$.  
  Therefore $N$ is the only minimal normal subgroup of $G$, $N\nleq \Phi 
(G) $ and $N$ is a $\pi$-group. Therefore $N\leq H$ and  $N=C_{G}(N)$ by 
\cite[A,  15.2]{DH}.  Since  $H\leq Z(N_{G}(H))$ and $H$ is 
a Hall  $\pi$-subgroup of $G$, $N=H$. Therefore $N\leq Z(G)$, which 
implies that $N=H=G$. This contradiction completes the proof of the lemma.

\section{Proof of Theorem A}

\

 Recall that $G$ is a $PST$-group if and only if  $G=D\rtimes M$, where $D=G^{\cal N
}$
 is abelian Hall subgroup of $G$ and  
every  element $x\in M$ induces a power automorphism on $D$ \cite{Agr}.   Therefore 
the implication
  (i) $\Rightarrow$ (ii) is a direct corollary of Theorem B.

    Now suppose that $G=D\rtimes M$, where $D=G^{\cal N
}$,  is a soluble 
$PST$-group.   Let $H$  be any  subgroup of $G$ and $S$ a Hall $\pi '$-subgroup of $H$.  
Since $G$ is soluble, we may 
assume without loss of generality  that $S\leq M$. Hence 
$H=(D\cap H)(M\cap H)=(D\cap H)S$ and $D\cap H$ is normal in $G$. 
 Let $\pi _{1}= \pi (S)$. Let $A$ be a Hall $\pi _{1}$-subgroup
 of $M$ and $E$ a complement to $A$ in  $M$.  Then $E\leq 
C_{G}(S)$. Therefore $G= DM=DAE=N_{G}(H)(DA)$ and every subgroup $L$ of $DA$ satisfying
 $(|H|, |L|)=1 $ is    contained in $D$.   Thus  $H$ 
  is quasipermutablein $G$. Thus (ii) $\Rightarrow$ (iii).

(iv) $\Rightarrow$ (ii) By Theorems C and D, $G$ is supersoluble and  $D$
 is  a Hall subgroup of $G$.
Therefore $G=D\rtimes W$, where $W$ is a Hall $\pi'$-subgroup of $G$. By hypothesis, $W$ is 
quasipermutable in $G$. Now arguing similarly as in the proof of Theorem 
B  one can show that $D$ is abelian and every subgroup of  $D$ is normal in $G$.  Therefore 
$G$ is a $PST$-group.

\section{Final remarks }

\

1. The subgroup $S_{3}$ is quasipermutable, $S$-propermutable and  not 
propermutable in    $S_{4}$. If $H$ is the subgroup of order 3 in $S_{3}$, 
then $H$ is $S$-quasipermutable and not quasipermutable in $S_{4}$.

2.  Arguing similarly to  the proof of Theorem   A  one can prove the following fact.

{\bf Theorem  5.1.}   {\sl Suppose that   $G$ is soluble  and let
 $\pi =\pi (G^{\cal N})$. Then $G$ is a $PST$-group if and only if 
 every subnormal $\pi$-subgroup and a Hall   $\pi'$-subgroup of $G$
 are  propermutable  in  $G$. }

3.    If $G$ is metanilpotent, that is $G/F(G)$ is nilpotent, 
then for every Hall subgroup $E$ of  $G$ we have $G=N_{G}(E)F(G)$.
Therefore, in this case, every characteristic
 subgroup of every Hall  subgroup of $G$ is $S$-propermutable in $G$.  In
particular, every Hall subgroup of every  supersoluble group is
$S$-propermutable. This observation makes natural the following
question: {\sl What is the structure of $G$ under the hypothesis
that every Hall subgroup of $G$ is propermutable in $G$ ?} Theorem B 
 gives an answer to this question.

4.  Every maximal subgroup of a supersoluble group is quasipermutable. Therefore, in fact,  
 Theorem  A   shows  that  the class of all soluble groups in which 
quaipermutability   is a  transitive relation   coincides with the 
class of all soluble $PST$-groups. 

5.      We say that $G$ is a \emph{$SQT$-group} if  $S$-quasipermutability     is 
a transitive relation  in $G$.   Arguing similarly to  the proof of 
Theorem   A  one can prove the following fact.  

{\bf Theorem 5.2.}  {\sl  A soluble group $G$ is an  $SQT$-group if and only if 
  $G=D\rtimes M $ is supersoluble,
 where $D$ and $M$
 are Hall nilpotent subgroups of $G$ and  the index
 $|G:DN_{G}(H\cap D)|$ is a 
$\pi (H)$-number  for every subgroup $H$ of $G$.  }

6.   A subgroup $H$ of $G$ is called  \emph{$SS$-quasinormal}    \cite{shirong}
 (\emph{semi-normal} \cite{Su}) in $G$ provided   $G$ has a subgroup $B$ 
such that $HB=G$ and  $H$ permutes with all  Sylow 
subgroups ($H$ permutes with all     subgroups, respectively)  of $B$.

It is clear that every $SS$-quasinormal subgroup is 
$S$-propermutable  and every  semi-normal   subgroup is propermutable.
Moreover, there are simple examples (consider, for example, the group
 $C_{7}\rtimes \text{Aut} (C_{7})$, where $C_{7}$ is a group of order 7) which show that,
 in general,  the class of all 
$S$-propermutable subgroups of $G$ is wider than the class of all its
$SS$-quasinormal subgroups and   the class of all 
propermutable subgroups of $G$ is wider than the class of all its
semi-normal  subgroups.  
Therefore Proposition  covers main results (Theorems 1.1--1.5) in  \cite{shirong}.

7.  Theorem 3.1 is used in the proof of Theorem B.
 From this result  we also get

{\bf Corollary 5.3}  (See \cite[Theorem 5.4]{8}). {\sl Let $H$ be 
  a Hall semi-normal  subgroup
of $G$.  If  $p
> q $ for all primes $p$ and $q$ such that $p$ divides $|H|$ and $q$
divides $|G:H|$,  then $H$ is normal in $G$.     }

{\bf Corollary 5.4}  (See \cite[Theorem]{GuoS}).  {\sl  Let $P$ be a
 Sylow $p$-subgroup of $G$. If $P$ is semi-normal in $G$, then the
following statements hold: }

(i) {\sl $G$ is $p$-soluble and $P'\leq O_{p}(G)$.}

(ii) {\sl   $l_{p}(G) \leq 2$. }

(iiii {\sl If for some prime $q\in  \pi'$ a Hall $p'$-subgroup
 of $G$
 is $q$-supersoluble, then $G$ is   $q$-supersoluble. }

{\bf Corollary 5.5}  (See  \cite[Theorem 3]{podg}).  {\sl   If a
Sylow $p$-subgroup $P$ of
 $G$,  where $p$ is the largest prime dividing $|G|$,
is semi-normal   in $G$, then $P$ is normal in $G$.}

\end{document}